\documentclass[12pt]{amsart}
\usepackage{fancyhdr}
\usepackage{amssymb,latexsym,amscd}
\usepackage[cp850]{inputenc}
\newcommand{\ds}{\displaystyle}
\newcommand{\cU}{\mathcal{U}}

\newcommand{\cE}{\mathcal{E}}
\newcommand{\cF}{\mathcal{F}}
\newcommand{\cM}{\mathcal{M}}
\newcommand{\lra}{\longrightarrow}
\newcommand{\hra}{\hookrightarrow}
\newcommand{\ra}{\rightarrow}
\newcommand{\lms}{\longmapsto}

\newcommand{\vf}{\varphi}

\newcommand{\gt}{\theta}
\newcommand{\ga}{\alpha}
\newcommand{\gb}{\beta}
\newcommand{\gc}{\gamma}

\newcommand{\cL}{\mathcal L}
\newcommand{\CC}{\mathbb C}

\newcommand{\ZZ}{\mathbb Z}

\newcommand{\PP}{\mathbb P}

\theoremstyle{plain}
\newtheorem{thm}{Theorem}[section]
\newtheorem{lem}[thm]{Lemma}
\newtheorem{prop}[thm]{Proposition}

\newtheorem{rem}[thm]{Remark}

\begin{document}
\title[Maximal Subbundles]{Maximal Subbundles and Gromov-Witten Invariants}
\author{H. Lange}
\author{P.E. Newstead}
\dedicatory{To C. S. Seshadri on his 70th birthday}
\address{Mathematisches Institut\\
              Universit\"at Erlangen-N\"urnberg\\
              Bismarckstra\ss e $1\frac{ 1}{2}$\\
              D-$91054$ Erlangen\\
              Germany}
\email{lange@mi.uni-erlangen.de}
\address{Department of Mathematical Sciences\\University of Liverpool\\
           Peach Street\\Liverpool L69 7ZL\\
           U.K.}
\email{newstead@liv.ac.uk}
\thanks{Supported by DFG Contracts Ba 423/8-1. Both authors are members
of the research group VBAC (Vector Bundles on Algebraic Curves), which is
partially supported by EAGER (EC FP5 Contract no. HPRN-CT-2000-00099) and
by EDGE (EC FP5 Contract no. HPRN-CT-2000-00101).}
\begin{abstract}
Let $C$ be a nonsingular irreducible projective curve
of genus $g\ge2$ defined over
the complex numbers. Suppose that $1\le
n'\le n-1$ and $n'd-nd'=n'(n-n')(g-1)$. It is known that, for the
general vector bundle
$E$ of rank $n$ and degree $d$, the maximal degree of a subbundle
of $E$ of rank $n'$ is $d'$ and that there are finitely many such
subbundles. We obtain a formula for the number of these maximal
subbundles when
$(n',d')=1$. For $g=2$, $n'=2$, we evaluate this formula explicitly.
The numbers computed here are Gromov-Witten invariants in the sense
of a recent paper of Ch. Okonek and A. Teleman (to appear in Commun.
Math. Phys.) and our results answer a question raised in
that paper.
\end{abstract}
\maketitle

%%%%%%%%%%%%%%%%%%%%%%%%%%%%%%%%%%%%%%%%%%%%%%%%%%%%%%%%%%

\section{Introduction}

Let $C$ be a non-singular irreducible projective curve of genus $g\ge2$, defined
over the complex numbers, and let
$E$ be a vector bundle over
$C$ of rank
$n$ and degree $d$. In describing the structure of $E$, it is
important to consider the set of subbundles $E'$ of $E$ of any fixed
rank $n'$ and degree $d'$. In particular, for any $n'$, there is a
maximum value of $d'$ for which such subbundles exist; the
corresponding subbundles are called {\it maximal subbundles} of $E$.
It is convenient here to write
$$s(E,E')=n'd-nd'.$$ We can then put
$$s_{n'}(E)=\min s(E,E'),$$the minimum being taken over all
subbundles $E'$ of rank $n'$. A subbundle $E'$ is then a maximal
subbundle if and only if $s(E,E')=s_{n'}(E)$.

Note in particular that the bundle $E$ is stable if and only if
$s_{n'}(E)>0$ for all $n'$, $1\le n'\le n-1$. In general we call the
$s_{n'}(E)$ {\it degrees of stability}  and note that, for each $n'$,
we can define a stratification of the moduli space ${\mathcal M}(n,d)$
of stable bundles of rank $n$ and degree $d$ by the locally closed
subsets
$$U_{n',s}(n,d)=\{E\in{\mathcal M}(n,d)|s_{n'}(E)=s\}.$$
Note that $U_{n',s}(n,d)$ can be non-empty only if $s>0$ and $s\equiv
n'd\bmod n$.

It was proved in \cite{11} that $s_{n'}(E)\le n'(n-n')g$ and in fact it is
known \cite[Th\'eor\`eme 4.4]{6} that the maximum value of $s$ that occurs is
given by
$$s=n'(n-n')(g-1) + \varepsilon,$$where $\varepsilon$ is the unique
integer $0\le\varepsilon\le n-1$ such that $s\equiv n'd\bmod n$. The
corresponding $U_{n',s}(n,d)$ is the unique open stratum of the
stratification.

When $s=n'(n-n')(g-1)$, a dimension counting argument shows that one
can expect a general $E\in {\mathcal M}(n,d)$ to have a finite number
$m_{n'}(E)$ of maximal subbundles of rank $n'$, the number being
independent of $E$ provided that it is chosen sufficiently generally.
In the simplest case, when $n'=1$, the number is known; in
fact,$$m_1(E)=n^g$$(see \cite{5,9,17} for $n=2$, \cite{15} and
\cite[Proposition
3.9]{14} in general).

Our object in this paper is to obtain a formula for $m_{n'}(E)$ in
the case $(n',d')=1$. Indeed we use Grothendieck-Riemann-Roch and the
Porteous formula to express ${\mathcal M}(n',d')$ as the top Chern
class of a certain virtual bundle (Theorem 3.1). When
$g=2$ and
$n'=2$, we evaluate this Chern class and give an explicit formula for
$m_2(E)$ (Theorem 4.1).

Our result can be interpreted in terms of Gromov-Witten theory .
Indeed, in \cite{14}, the formula for $m_{n'}(E)$ is obtained as an
application of a theorem on Gromov-Witten invariants, and the
question of computing these invariants when $n'>1$ is raised. Our
results give an answer to this question. Our methods, however, are
similar to those of \cite{15}. Notice that these invariants were also defined
in the last section of \cite{3}.

It should also be noted that there is a mistake in \cite{16}, where it is
claimed that, when $s\le n'(n-n')(g-1)$, the general bundle in
$U_{n',s}(n,d)$ has only one maximal subbundle of rank $n'$. For
$s=n'(n-n')(g-1)$, this is not correct, and indeed the proof given in
\cite{16} fails also for $s<n'(n-n')(g-1)$. However the result itself is
true in this case; we provide a proof of this in Theorem 2.3. It
should be emphasised that this does not invalidate any other result
in \cite{16}.

The research for this paper was carried out during a visit by the
second author to the Mathematisches Institut der Universit\"at
Erlangen-N\"urnberg in October/November 2001. He is grateful to the
DFG for funding this visit and to the Institut for its hospitality.

%%%%%%%%%%%%%%%%%%%%%%%%%%%%%%%%%%%%%%%%%%%%%%%%%%%%%%%%%%%

\section{Uniqueness of Maximal Subbundles}

According to \cite[Theorem 0.1]{16} $U_{n',s}(n,d)$
is locally closed in $\mathcal M (n,d)$ and irreducible of dimension
$$
\dim U_{n',s}(n,d) = (n^2 - n'(n - n'))(g - 1) + s + 1.
$$
Hence it makes sense to speak of a general vector bundle in $U_{n',s}(n,d)$.\\
\begin{lem}
Let $E$ be a general bundle in $U_{n',s}(n,d)$ with $ s \leq n'(n - n')(g - 1)$.
Then $E$ possesses only finitely many maximal subbundles $E'$.
For any one of them both
$E'$ and $E'' = E/E'$ are general (as points of the moduli spaces $\mathcal M
(n',d')$ and $\mathcal M(n-n',d-d')$ respectively).
\end{lem}

\noindent
\begin{rem}
{\rm Lemma 2.1 has been proved in \cite[Claim p. 495 and Lemma 1.4]{16}. Similar
statements have been given in
\cite{4,6,8}. We include the proof for the convenience of the reader.}
\end{rem}

\begin{proof} We show first that $E'$ and $E''$ are both general. If not, then either
$E'$ or $E''$ (or both) depends on fewer parameters than the dimension
of the corresponding moduli space; this is true irrespective of whether $E'$ and $E''$
are stable. Note also that $h^0(E''{}^*\otimes E')=0$ since $E$ is stable.
So the nontrivial extensions
$$
0 \lra E' \lra E \lra E'' \lra 0
$$
depend on a number $\nu$ of parameters with
\begin{align*}
\nu &< \dim \mathcal M(n',d') + \dim \mathcal M(n-n',d-d') +
h^1(E''{}^* \otimes E') -1\\
&= (g-1)(n'{}^2 + (n-n')^2 + n'(n-n')) + s + 1\\
&= (g-1)(n^2 - n'(n-n')) + s + 1\\
&= \dim U_{n',s}(n,d).
\end{align*}
This contradicts the generality of $E$.

If the general $E$ in $U_{n',s}(n,d)$ possesses infinitely many maximal
subbundles, then the same inequality holds, again a contradiction.
\end{proof}

If $s < n' (n-n')(g-1)$, we can be more precise:\\

\begin{thm}
A general vector bundle $E \in U_{n',s} (n,d)$ with $s < n' (n-n')
(g-1)$ admits only one maximal subbundle of rank $n'$.
\end{thm}

{\it Proof.} For convenience we write $n''=n-n'$, $d''=d-d'$.
Passing if necessary to the dual bundle, we may assume that $n'\le
n''$. Moreover assume first, slightly more generally, that $s \le
n'(n-n')(g-1)$, since we want to deduce also something for $s
= n'(n-n')(g-1)$.

Fix a maximal subbundle $E' \in \cM(n',d')$ of $E$ and write $E''
= E/E' \in \cM (n'',d'')$. So
$$
s=n'd'' - n''d'.
$$
Suppose $F$ is a second maximal subbundle of $E$. The composed map
$\vf:\ F\hra E \ra E''$ is nonzero. Let $G\subset E''$ be its
(sheaf theoretic) image.
Then we have the following commutative diagram
$$
\begin{array}{ccccccccc}
0 & \ra & H & \ra & F & \ra & G & \ra & 0 \\
& & \downarrow_j & & \downarrow & & || & & \\
0 & \ra & E' & \ra & \widetilde{E} & \ra & G & \ra & 0 \\
& & || & & \downarrow & & \downarrow_i & & \\
0 & \ra & E' & \ra & E & \ra & E'' & \ra & 0 \\
\end{array}
$$
where the middle exact sequence is the push out of the upper exact
sequence and the pull back of the lower sequence and the composition
$F\lra E''$ is the map $\vf$. Note that $i$ and $j$ are injective
(as maps of sheaves). Conversely, given such a diagram with $rk(F) = n'$ and $deg
(F) =d'$, then $F$ is a second maximal subbundle of $E$ of rank $n'$.
Let $n_G$ and $d_G$ denote the rank and degree of $G$, $1\le n_G\le n'$.

Let $A_{n_G,d_G} (E'')$ denote the set of subsheaves of $E''$ of rank
$n_G$ and degree $d_G$, and define similarly $A_{n_H,d_H} (E')$.
Note that $A_{n_G,d_G} (E'')$ and $A_{n_H,d_H} (E')$ can be given a
scheme structure by identifying them with Quot schemes. Finally
denote by $i^*$ the pull back map $H^1 (E''{}^* \otimes E') \ra H^1
(G^*\otimes E')$. According to the above remarks, it suffices to show
that, if $s< n'(n-n')(g-1)$, then
$$
\dim A_{n_G,d_G} (E'') + \dim A_{n_H,d_H} (E') + h^1 (G^* \otimes H)
+ \dim Ker (i^*) < h^1 (E''{}^* \otimes E') \eqno (1)
$$
But, by Lemma 2.1 and \cite[Theorem 0.2]{16},
$$
\begin{array}{rcl}
\dim A_{n_G,d_G} (E'') &=& n_G d'' - n'' d_G - n_G (n''-n_G) (g-1) \\
\dim A_{n_H,d_H} (E') &=& n'd_G - n_G d' -(n'-n_G) n_G (g-1). \\
\end{array}
$$
Moreover, since $i^*$ is surjective and $h^0 (G^*\otimes H)=0,\ F$
being stable,
$$
\begin{array}{rcl}
h^1 (G^* \otimes H) &=& n' d_G - n_G d' + n_G (n' -n_G) (g-1) \\
\dim Ker (i^*) &\le& n'(n''-n_G) (g-1) + s + n_G d' - n'd_G. \\
\end{array}
$$
Hence inequality (1) would follow from
$$
(n'n'' + n_G (n_G - n'-n''))(g-1)+n_G d'' - n'' d_G + n' d_G - n_G
d' + s < n'n''(g-1) +s
$$
which is equivalent to
$$
n_G (n_G -n'-n'')(g-1) + n_G (d''-d') + (n'-n'')d_G < 0. \eqno (2)
$$
But Lemma 2.1 implies that $n_G (n'-n_G)(g-1) \le s_{n'-n_G}(F) \le n'd_G -
n_Gd'$ which gives
$$
d_G \ge \frac{n_G}{n'} \left[(n'-n_G)(g-1)+d' \right].
$$
Hence the left hand side of (2) is less than or equal to (here we use
the fact that $n'\le n''$)
%\begin{align*}
$$
\begin{array}{rlr}
\ds{\frac{n_G}{n'}} & \ds{\left[n''(n_G-2n')(g-1) + d'' n'-n''d' \right]}&\\
&\ds{= \frac{n_G}{n'} \left[n''(n_G-n')(g-1) + s - n'n''(g-1)\right]}& \\[.2cm]
&\ds{\le \frac{n_G}{n'} \left[s-n'n''(g-1)\right]\ (\hbox{since }\ n_G
\le n')}& \\[.2cm]
&\ds{< 0\ \hbox{ for }\ s < n'n'' (g-1).}& \\
\end{array}
%\end{align*}
$$
$$
\eqno \square
$$
\noindent
Suppose now $s=n'(n-n')(g-1)$. According to Lemma 2.1 there are only
finitely many maximal subbundles for general $E$. We obtain as a
consequence of the proof of Theorem 2.3
\begin{prop}
Let $E$ be a general vector bundle on $C$ with $s_{n'}(E) = n'(n-n')(g-1)$
and $n'\le (n-n')$. If $E'$ and $F$ are $2$ maximal subbundles of rank
$n'$ of $E$, the composed map $F\ra E \ra E/E'$ is injective.
\end{prop}
\begin{proof}
For the proof just note that in the proof of Theorem 2.3 all
inequalities have to be equalities. In particular we must have
$n_G = n'$, which implies the assertion.
\end{proof}

%%%%%%%%%%%%%%%%%%%%%%%%%%%%%%%%%%%%%%%%%%%%%%%%%%%%%%%%%

\section{Maximal subbundles of general bundles}

Let $E$ be a general stable bundle of rank $n$ and degree $d$
on $C$ and let $n'$ be an integer, $1\le n' \le n-1$, such that
$$
s_{n'} (E) = n' (n-n')(g-1) = n'd-nd'.
$$
We know by Lemma 2.1 that $E$ has only a finite number of maximal
subbundles. Our object in this section is to obtain a formula
for this number when $(n',d')=1$.

We denote by $M_{n'} (E)$ the set of maximal subbundles of $E$
of rank $n'$. This set can be given a scheme structure by identifying
it with Grothendieck's scheme $Quot_{n-n',d''} (E)$
of quotients of $E$ of rank $n-n'$ and degree
$$
d''= d-d' = \frac{1}{n} \left( (n-n')d+s_{n'} (E) \right)
$$

Now let $J$ denote the Jacobian of $C$ consisting of line bundles
of degree $0$ on $C$, and let $M_0$ denote the moduli space of
stable bundles of rank $n'$ and with fixed determinant of degree
$d'$ on $C$. The map
$$
\pi:\ J\times M_0 \lra {\mathcal M}(n',d') :\ (L,F) \lms L\otimes F
$$
is an unramified covering of degree $n'{}^{2g}$.

Let $\cL$ denote a Poincar\'e bundle on $C\times J$ and $\cU$ a
universal bundle on $C\times M_0$ (recall that we are
assuming $(n',d')=1$, so $\cU$ exists). Let $K$ denote the
canonical bundle on $C$ and $q:\ C\times J\times M_0 \lra J\times
M_0$ and $p:\ C\times J \times M_0 \lra C$ be the canonical
projections. We denote the pullbacks of $\cU$ and $\cL$ to $C\times
J\times M_0$ by the same symbols and the restrictions of $\cU$ and
$\cL$ to $\{x\}\times J\times M_0$ by $\cU_x$ and $\cL_x$.\\

\begin{thm}
Let $E$ be a general bundle of rank $n$ and degree $d$, and $n'$
an integer with $1\le n' \le n-1$. Suppose that $n'd-nd' = n'
(n-n')(g-1)$ and that $(n',d')=1$. Then the number of maximal
subbundles of $E$ of rank $n'$ is given by
$$
 m_{n'} (E) = \frac{1}{n'{}^{2g}} c_{\rm top} \left(
\cF - \cE\right) [J\times M_0]
$$
where $\cE = q_* (\cU^*\otimes \cL^* \otimes p^* (E\otimes K))$ and
$\cF=\bigoplus_{x\in D}(\cU_x^*\otimes\cL_x^*)\otimes\CC^n$, where $D$ is a
smooth canonical divisor on $C$.
\end{thm}

\noindent
{\it Proof.}
Note first that the conclusion is equivalent to
$$
\# \pi^{-1} (M_{n'} (E)) = c_{\rm top} \left(
\cF - \cE\right) [J\times M_0].
$$
Consider the exact sequence
$$
0 \lra \cU^*\otimes \cL^* \otimes p^*E \stackrel{\otimes D}{\lra}
\cU^*\otimes \cL^* \otimes p^* (E \otimes K) \lra\cU^*\otimes \cL^*
|_{p^*(D)}\otimes
\CC^{n} \lra 0
$$
where $D$ denotes any smooth canonical divisor. Note that for any
maximal subbundle $E'$ of $E$,
$$
\deg (E'{}^* \otimes E \otimes K) = n'd - nd' + nn' (2g-2) >
nn' (2g-2).
$$
Since any such $E'$ is stable by Lemma 2.1, $E'{}^* \otimes E
\otimes K$ is semistable, so
$$
H^1 (E'{}^* \otimes E\otimes K) =0.
$$
Hence $R^1 q_* (\cU^*\otimes \cL^* \otimes p^* (E\otimes K)) =0$
and we have an exact sequence
$$
0 \lra q_* (\cU^*\otimes \cL^* \otimes p^* E) \lra \cE
\stackrel{\vf}{\lra} \cF \lra R^1 q_* (\cU^*\otimes \cL^* \otimes p^*E)
\lra 0
$$
with $\cE$, $\cF$ as in the statement of the theorem.
Note that, by Riemann-Roch,
$$
\begin{array}{rcl}
rk\ \cE &=& n'd-nd' +nn' (g-1) \\
&=& 2nn' (g-1) -n'{}^2 (g-1) \\
&=& rk\ \cF - n'{}^2 (g-1). \\
\end{array}
$$
\\

\noindent
{\bf Claim:}
$\pi^{-1} (M_{n'} (E)) =\{ (L,F) \in J\times M_0| rk\ \vf_{(L,F)}
\le rk\ \cE -1\}$.

\medskip
The claim is a consequence of the following lemma and Lemma 2.1.

\begin{lem}
Let $E$ be a general bundle. Suppose $n'd-nd'=n'(n-n')(g-1)$ and let
$E'$ be a stable bundle of rank $n'$ and degree $d'$. Then any
nonzero homomorphism $E' \lra E$ is an injection onto a
subbundle of $E$. Moreover $h^0 (E'{}^* \otimes E) =1$.
\end{lem}

\begin{proof}
If $\psi:\ E'\lra E$ is injective then  $Im \psi$ must be saturated
since it has the same degree as a maximal subbundle of $E$.
Suppose that $\psi$ is not injective and let $E_1 =Im \psi$. Since
$E'$ is stable, we have
$$
\frac{d_1}{n_1} > \frac{d'}{n'} = \frac{d}{n} - \frac{n-n'}{n}
(g-1) > \frac{d}{n} - \frac{n-n_1}{n} (g-1),
$$
contradicting the generality of $E$. Finally, if $h^0 (E'{}^* \otimes
E) \ge 2$, the morphism
$$\PP (H^0(E'{}^*\otimes E)) \lra M_{n'} (E)$$
is non-constant, since $E'$ is stable. So $M_{n'} (E)$ is infinite,
another contradiction.
\end{proof}

This establishes the claim. It follows by the Porteous formula (see
\cite[p.86]{1}) that the fundamental class of $\pi^{-1} (M_{n'}
(E))$ in $H^* (J\times M_0,\ZZ)$ is given by
$$
[\pi^{-1} (M_{n'} (E))] = c_{\rm top} (\cF - \cE)[J\times M_0].
$$
To complete the proof of the theorem, we must show that $M_{n'}
(E)$ is smooth.

Given the identification of $M_{n'} (E)$ with $Quot_{n-n',d''}
(E)$, it is sufficient to prove

\begin{lem}
Let $E$ be a general bundle and $E'$ any subbundle of $E$. Then
$$
H^1 (E'{}^* \otimes E/E') =0.
$$
\end{lem}

\begin{proof}
(due essentially to Laumon \cite {10}, see also \cite{15}). If
$H^1 (E'{}^* \otimes E/E') \ne 0$, then by Serre duality there exists
a nonzero homomorphism $E/E' \lra E' \otimes K$. We thus have a non-zero
homomorphism
$$
E \lra E/E' \lra E' \otimes K \lra E\otimes K.
$$
This homomorphism is clearly nilpotent. But, since $E$ is general,
it is very stable (see \cite[Proposition 3.5]{10}), i.e.\ it
admits no nilpotent homomorphism $E\lra E\otimes K$ different from
zero. This contradiction establishes the lemma and Theorem 3.1.
\end{proof}

\begin{rem}
{\rm It may be of interest to identify the precise generality
conditions required on $E$. In fact we need
\begin{enumerate}
\item[(i)] $s_{n'}(E)=n'(n-n')(g-1)$;
\item[(ii)] $s_{n_1}(E)\ge n_1(n-n_1)(g-1)$ for $1\le n_1\le n'-1$;
\item[(iii)] every maximal subbundle of $E$ of rank $n'$ is stable;
\item[(iv)] $M_{n'}(E)$ is finite;
\item[(v)] $E$ is very stable.
\end{enumerate}
The application of the Porteous formula, however, requires only the
assumption that the relevant degeneracy locus of $\vf$ is
finite and it then gives the length of this locus in its natural
scheme structure. In particular, if we assume only conditions (i),
(ii) and (iv), the formula counts (with the appropriate
multiplicities) the number of stable maximal subbundles of $E$ of
rank $n'$.}
\end{rem}
\begin{rem}
{\rm Since $(n',d')=1$ and $n'd-nd'=n'(n-n')(g-1)$, it follows that
$n'$ divides $n$ and we can then solve for $d$ in terms of $n,n',d'$.
Theorem 3.1 does therefore give an answer to the question of
\cite[section 3.4, Remark 4]{14} whenever $(n',d')=1$ (or, in the notation
of \cite{14}, $(r,d)=1$). In \cite{14} there is no smoothness assumption on the
Quot scheme, so we require only conditions (i)--(iv) of Remark 3.4;
the formula of Theorem 3.1 then gives the length of the finite scheme
$Quot_{n-n',d''}(E)$.}
\end{rem}

\begin{rem}
{\rm To apply this theorem, we note that the Chern classes of $\cU$ and
$\cL$ are well known (see \cite{2} for $\cU$ and \cite{1} for $\cL$).
From this the
Chern character of $\cU^* \otimes \cL^* \otimes p^* (E\otimes K)$ and
hence, using Grothendieck-Riemann-Roch, of $\cE$ can be computed.
The Chern character of $\cF$ is also easily computed, hence  that
of $\cF-\cE$.
Converting this to Chern classes, we obtain $c_{\rm top}
(\cF-\cE)$. We now need to know only the intersection numbers
of $J\times M_0$. In principle, these can be calculated. In
practice, it may be a difficult computation. In the final
section of this paper, we will carry out the computation in the
simplest case not previously known.}
\end{rem}
One could of course have stated the theorem using the cohomology
of $M(n',d')$ rather than that of $J\times M_0$. We have chosen the
latter, because (at least in the case $n'=2$) it is easier to work
with.

%%%%%%%%%%%%%%%%%%%%%%%%%%%%%%%%%%%%%%%%%%%%%%%%%%%%%%%%%

\section{The case $g=2,\ n'=2$}

\begin{thm}
Let $E$ be a general bundle of rank $n$ and degree $d$ on a
curve $C$ of genus $2$. Suppose that $n\ge 4$ is even and
$2d+4 \equiv 0 \mod n$ with $\frac{2d+4}{n}$ odd. Then
$$
 m_2 (E) =\frac{n^3}{48} (n^2 + 2).
$$
\end{thm}

\begin{rem}
{\rm Note that the right hand side of this expression is always a positive
integer for $n$ even. The numerical conditions are precisely
what is needed to allow $d'$ to be odd. When they are satisfied
we can take $d'$ to be any odd integer and then}
$$
d =\frac{n}{2} d' + n-2.
$$
\end{rem}
\vspace{0.5cm}

Before starting on the proof, we recall some properties of the
cohomology ring $H^* (C\times J \times M_0,\ZZ)$. The abelian group
$H^1 (C,\ZZ)$ has a natural symplectic structure. Let $f$ denote the
positive generator of $H^2 (C,\ZZ)$. We can then choose a
symplectic basis $e_1,e_2,e_3,e_4$ for $H^1 (C,\ZZ)$ such that
$$
e_1e_3 = e_2 e_4 = -f \eqno (1)
$$
and all the other products $e_ie_j$ for $i<j$ are 0. On $C\times J$
we can normalise $\cL$ so that
$$
c(\cL) = 1 + \xi_1 \eqno (2)
$$
where $\xi_1 \in H^1 (C,\ZZ) \otimes H^1 (J,\ZZ) \subset H^2 (C\times
J)$ can be written as
$$
\xi_1 = \sum_{i=1}^{4} e_i \otimes \vf_i \eqno (3)
$$
with $\vf_i \in H^1 (J,\ZZ)$. We have (see \cite[p.335]{1})
$$
\xi_1^2 = -2 \gt f,\quad \gt^2 [J] =2. \eqno (4)
$$

The variety $M_0$ is given by the intersection of two quadrics
in $\PP_5$ (see \cite{12}). From this we can see that its cohomology
groups are
$$
\ZZ\, ,\ 0 \, ,\ \ZZ\, ,\ \ZZ^4\, ,\ \ZZ\, ,\ 0 \, ,\ \ZZ.
$$
If $\ga$ is the positive generator of $H^2 (M_0, \ZZ)$, then $\ga^3
[M_0] =4$, so by (4)
$$
\ga^3 \gt^2 [J\times M_0] =8. \eqno (5)
$$

Taking $d'=1$ and normalising as in \cite{13} we have
$$
c_1 (\cU) = \ga + f,\quad c_2 (\cU) = \chi + \xi_2 + \ga f\eqno (6)
$$
with $\chi \in H^4 (M_0,\ZZ)$ and $\xi_2 = H^1 (C,\ZZ) \otimes
H^3 (M_0,\ZZ)$. In the notation of \cite{13} and \cite{7}, $4\chi
= \ga^2 - \gb$ and we write
$$
\xi_2^2 = \gc f,\quad \gc \in H^6 (M_0,\ZZ). \eqno (7)
$$
Moreover, when $g=2$, the relations of the Theorem of \cite{7}
give
$$
\ga^2 + \gb =0,\quad \ga^3 + 5 \ga\gb + 4 \gc =0.
$$
Hence $\gb = -\ga^2 \, ,\ \gc =\ga^3$, and we have from (6)
$$
c_2 (\cU) = \frac{\ga^2}{2} + \xi_2 + \ga f,\quad \xi_2^2
= \ga^3 f. \eqno (8)
$$

We define also $\Lambda \in H^1 (J,\ZZ) \otimes H^3 (M_0,\ZZ)$ by
$$
\xi_1\xi_2 = \Lambda f. \eqno (9)
$$
In $H^* (C\times J \times M_0,\ZZ)$, we have for dimensional reasons
(noting that $H^5 (M_0,\ZZ)=0$) that the following classes are
all zero:
$$
f^2\, ,\ \xi_1^3\, ,\ \alpha^4 \, ,\ \xi_1 f \, ,\ \xi_2 f\, ,\ \ga\xi_2\, ,\
\ga \Lambda\, ,\ \gt^2 \Lambda. \eqno (10)
$$

Finally we need to compute $\gt \Lambda^2$. We can write
$$
\xi_2 = \sum_{i=1}^4 e_i \otimes \psi_i \quad \hbox{ with }\
\psi_i \in H^3 (M_0,\ZZ).
$$
Using this, (1) and (3), we have
$$
\xi_1 \xi_2 = (\vf_1 \psi_3 + \vf_2 \psi_4 - \vf_3 \psi_1 - \vf_4
\psi_2) f.
$$
So by (9)
$$
\Lambda = \vf_1 \psi_3 + \vf_2 \psi_4 - \vf_3 \psi_1 -
\vf_4 \psi_2.
$$
Also from (3) and (4),
$$
\gt = -\vf_1\vf_3 - \vf_2\vf_4.
$$
Hence $\gt^2 = 2\vf_1\vf_3\vf_2\vf_4$, so by (4)
$$
\vf_1\vf_3\vf_2\vf_4 [J] =1.
$$
Now
$$
\begin{array}{rcl}
\Lambda^2 &=& -2\vf_1\vf_2\psi_3\psi_4 + 2\vf_1\vf_3\psi_3\psi_1
+ 2 \vf_1\vf_4\psi_3\psi_2 \\
&+& 2 \vf_2\vf_3\psi_4\psi_1 + 2\vf_2\vf_4\psi_4\psi_2 -
2 \vf_3\vf_4\psi_1\psi_2 \\
\end{array}
$$
and
$$
\begin{array}{rcl}
\gt\Lambda^2 &=& -2\vf_1\vf_3\vf_2\vf_4 \psi_4\psi_2
- 2 \vf_2\vf_4\vf_1\vf_3\psi_3\psi_1 \\
&=& \gt^2 (\psi_1\psi_3 + \psi_2 \psi_4) \\
&=& \frac{1}{2} \gt^2 \gc \qquad \qquad (\hbox{by}\ (7))\\[.2cm]
&=& \frac{1}{2} \gt^2 \ga^3 \qquad \qquad (\hbox{since }\
\gc=\ga^3). \\
\end{array}
$$
So by (5)
$$
\gt\Lambda^2 [J\times M_0] = 4. \eqno (11)
$$

\begin{proof}[Proof of Theorem 4.1.] Tensoring by a line bundle over $X$
if necessary, we can suppose $d'=1$.
By (2) and (4) we have
$$
ch (\cL) =1 + \xi_1 - \gt f
$$
while by (6), (8) and (10)
$$
ch(\cU) = 2 + (\ga+f) +(-\xi_2) +\left(-\frac{1}{12} \ga^3 -
\frac{1}{4} \ga^2 f\right).
$$
Moreover
$$
ch (p^* (E\otimes K)) \cdot p^* td (C) = n + \left( \frac{5}{2}
n-2\right) f.
$$
So
\begin{align*}
ch (\cU^* \otimes \cL^* \otimes p^* (E\otimes K)) & \cdot p^* td (C) =
2n + [(4n-4) f - n\ga - 2n\xi_1] \\
&+ \left[\left(-\left(\frac{5}{2} n-2\right)\ga - 2n\gt\right) f
+ n\ga\xi_1 - n\xi_2 \right] \\
&+\left[\left(\frac{n}{4} \ga^2 + n\Lambda + n\ga \gt\right) f +
\frac{n}{12} \ga^3 \right] \\
&+\left[\left(\frac{5}{24} n-\frac{1}{6}\right) \ga^3 f -
\frac{n}{12} \ga^3 \xi_1\right] + \left[ - \frac{n}{12} \ga^3
\gt f \right]. \\
\end{align*}
Hence by Grothendieck-Riemann-Roch
\begin{align*}
ch (\cE) = 4n-4 &+
\left[-\left(\frac{5}{2} n- 2\right)\ga - 2n \gt\right] +
\left[\frac{n}{4} \ga^2 +n\Lambda +n\ga \gt\right] \\
& + \left[\left(\frac{5}{24}n-\frac{1}{6} \right) \ga^3  \right] +
\left[ - \frac{n}{12} \ga^3 \gt \right]. \\
\end{align*}
From (6) and (8) we obtain easily
$$ch(\cF)=4n-2n\alpha+\frac{n}{6}\alpha^3.$$ So
\begin{align*}
ch (\cF-\cE) = 4 &+
\left[\left(\frac{1}{2} n- 2\right)\ga + 2n \gt\right] +
\left[-\frac{n}{4} \ga^2 -n\Lambda -n\ga \gt\right] \\
& + \left[\left(\frac{1}{6} -\frac{1}{24} n\right) \ga^3  \right] +
\left[\frac{n}{12} \ga^3 \gt \right]. \\
\end{align*}
Converting this into Chern classes, we obtain
$$
c_{\rm top} (\cF -\cE) = c_5
(\cF -\cE) = \left(\frac{1}{24} n^5 -
\frac{5}{12} n^3\right) \ga^3 \gt^2 + n^3 \gt \Lambda^2.
$$
By (5) and (11), we get
$$
c_{\rm top} (\cF -\cE) [J\times M_0]
= \frac{1}{3} n^5  + \frac{2}{3} n^3
$$
and Theorem 4.1 follows from Theorem 3.1.
\end{proof}

%%%%%%%%%%%%%%%%%%%%%%%%%%%%%%%%%%%%%%%%%%%%%%%%%%%%%%%%%%
\vspace{1cm}
%%%%%%%%%%%%%%%%%%%%%%%%%%%%%%%%%%%%%%%%%%%%%%%%%%%%%%%%%%

\end{document}